\documentclass[letterpaper,11pt]{article}
\usepackage{amsfonts}

\newtheorem{theo}{Theorem}%[section]

\newtheorem{lemma}{Lemma}[section]
\newtheorem{definition}[lemma]{Definition}

\newtheorem{claim}[lemma]{Claim}

\newcommand{\ignore}[1]{}

\newcommand{\qed}{\hspace*{\fill} \rule{7pt}{7pt}}

\begin{document}
\title{The Effect of Induced Subgraphs on Quasi-Randomness\thanks{A
preliminary version of this paper appeared in the Proc. of the $19^{th}$
Annual ACM-SIAM Symposium on Discrete Algorithms (SODA), ACM Press (2008),
789-798.}} \author{
Asaf Shapira
\thanks{School of Mathematics and College of Computing, Georgia Institute of Technology, Atlanta, GA 30332}
\and
Raphael Yuster
\thanks{Department of Mathematics, University of Haifa, Haifa
31905, Israel. E--mail: raphy@math.haifa.ac.il} }

\date{}
\maketitle

\setcounter{page}{1}
\begin{abstract}

One of the main questions that arise when studying random and
quasi-random structures is which properties ${\cal P}$ are such that
any object that satisfies ${\cal P}$ ``behaves'' like a truly random
one. In the context of graphs, Chung, Graham, and Wilson \cite{CGW}
call a graph {\em $p$-quasi-random} if it satisfies a long list of
the properties that hold in $G(n,p)$ with high probability, like
edge distribution, spectral gap, cut size, and more.

Our main result here is that the following holds for {\bf any}
fixed graph $H$: if the distribution of {\em induced} copies of
$H$ in a graph $G$ is close (in a well defined way) to the
distribution we would expect to have in $G(n,p)$, then $G$ is
either $p$-quasi-random or $\overline{p}$-quasi-random, where
$\overline{p}$ is the unique non-trivial solution of the
polynomial equation $x^{\delta}(1-x)^{1-\delta} =
p^{\delta}(1-p)^{1-\delta}$, with $\delta$ being the edge density
of $H$. We thus infer that having the correct distribution of
induced copies of any single graph $H$ is enough to guarantee
that a graph has the properties of a random one. The proof
techniques we develop here, which combine probabilistic, algebraic
and combinatorial tools, may be of independent interest to the
study of quasi-random structures.

\end{abstract}

\thispagestyle{empty}
\setcounter{page}{1}

\section{Introduction}\label{SecIntro}

\subsection{Background and basic definitions}

Quasi-random structures are those that possess
the properties we expect random objects to have with high
probability. The study of quasi-random structures is one of the most
interesting borderlines between discrete mathematics and theoretical
computer science, as they relate the problem of how to
deterministically construct a random-like object with the question
of when we can consider a single event to be a random one. Although
quasi-random structures have been implicitly studied for many
decades, they were first explicitly studied in the context of graphs
by Thomason \cite{T1,T2} and then followed by Chung, Graham, and
Wilson \cite{CGW}. Following the results on quasi-random graphs,
quasi-random properties were also studied in various other contexts
such as set systems \cite{CG}, tournaments \cite{CG1}, and
hypergraphs \cite{CG2}. There are also some very recent results on
quasi-random groups \cite{Go} and generalized quasi-random graphs
\cite{LS}. We briefly mention that the study of quasi-random
structures lies at the core of the recent proofs of Szemer\'edi's
Theorem \cite{Sztheo} that were recently obtained independently by
Gowers \cite{Go1,Go2} and by Nagle, R\"odl, Schacht and Skokan \cite{RS,NRS} and then
also by Tao \cite{T} and Ishigami \cite{Ish}. For more mathematical background on
quasi-randomness the reader is referred to the recent papers of
Gowers \cite{Go,Go1,Go2} and to the survey of Krivelevich and
Sudakov \cite{KS}.

Besides being intriguing questions on their own, results on
quasi-random objects also have applications in theoretical computer
science. The main point is that, while the classical definitions of
what it means for an object to be quasi-random are hard to verify,
some other properties, which can be proved to be equivalent, are
much easier to verify. The archetypal example of this phenomenon is
probably the spectral gap property of {\em expanders}. Expanders are
sparse graphs that behave like random sparse graphs in many aspects
(see \cite{HLW} for more details), and are one of the most widely
used structures in theoretical computer science. However, verifying
that a graph satisfies the classical definition of being an
expander, that is, that any cut has many edges, requires exponential
time. A very useful fact is that being an expander is equivalent to
the fact that the absolute value of the second eigenvalue of the adjacency matrix of the
graph is significantly smaller than the first eigenvalue (see also
Property ${\cal P}_3$ in Theorem \ref{theoCGW}). As eigenvalues can
be computed in polynomial time, this gives an efficient way to
verify that a sparse graph is an expander.

Throughout the paper, when we say that a graph $G$ has the ``correct''
number of copies of a graph $H$ as we would expect to have in $G(n,p)$,
we mean that the number of copies of $H$ in $G$ is $(1+o(1))n^p$, where as usual,
an $o(1)$ term represents an arbitrary function tending to $0$ as $n$ tends to infinity.
Let us now consider another example in which equivalence between different notions of quasi-randomness is useful,
this time on dense graphs. A natural notion of quasi-randomness
for a dense graph is that all subsets of vertices should contain
the ``correct'' number of edges as in $G(n,p)$. This property takes
exponential time to verify, but fortunately (see Theorem
\ref{theoCGW}), it turns out that this property is equivalent to the
property of having the ``correct'' number of edges and copies of the
cycle of length four in the entire graph!  As this property takes
only polynomial time to verify, this gives an efficient algorithm
for checking if a dense graph is quasi-random. This easily
verifiable condition was a key (implicit) ingredient in the work of
Alon, Duke, Lefmann, R\"odl and Yuster \cite{ADLRY}, who gave the first polynomial time
algorithm for Szemer\'edi's Regularity Lemma \cite{Sz}, whose
original proof was non-constructive.

Given the above discussion, one of the most natural questions that
arise when studying quasi-random objects, is which properties
``guarantee'' that an object behaves like a truly-random one. Our
main result in this paper establishes that for any single graph
$H$, if the distribution of the induced copies of $H$ in a graph
$G$ is ``close'', in some precise sense, to the one we expect to
have in $G(n,p)$, then $G$ is quasi-random. Previous studies
\cite{SSnind,SSind} of the effect of induced subgraph on
quasi-randomness that used a slightly weaker notion of
``closeness'', indicated that in some cases the distribution of
induced copies of a single graph $H$ is {\em not enough} to
guarantee that a graph is quasi-random. Therefore, the notion of
closeness that we use here is essentially optimal if one wants to
be able to deal with any $H$.

Before stating our main result we first discuss some previous
ones, which will put ours in the right context. The cornerstone
result on properties guaranteeing that a graph is quasi-random is
that of Chung, Graham, and Wilson \cite{CGW}, stated below, but
before stating it we need to introduce some notation. We will
denote by $e(G)$ the number of edges of a graph $G$. A {\em
labeled copy} of a graph $H$ in a graph $G$ is an injective
mapping $\phi$, from the vertices of $H$ to the vertices of $G$,
that maps edges to edges, that is $(i,j) \in E(H) \Rightarrow
(\phi(i),\phi(j)) \in E(G)$. So the expected number of labeled
copies of a graph $H$ in $G(n,p)$, is $p^{e(H)}n^h+o(n^h)$
where $h$ is the number of vertices of $H$
\footnote{Note that this is {\em not} the expected number of
unlabeled copies of $H$ in $G$, which is just the number of
labeled copies of $H$ divided by the number of automorphisms of
$H$. Therefore, all the results we mention here also hold when
considering unlabeled copies. We work with labeled copies
(induced or not) because we do not need to refer to the
automorphisms of $H$, and because it is easier to count labeled
copies than copies.}. A {\em labeled induced copy} of a graph $H$
in a graph $G$ is an injective mapping $\phi$, from the vertices
of $H$ to the vertices of $G$, that maps edges to edges, and non-edges
to non-edges, that is $(i,j) \in E(H) \Leftrightarrow
(\phi(i),\phi(j)) \in E(G)$. So the expected number of induced
labeled copies of a graph $H$ in $G(n,p)$, is
$\delta_H(p)n^h+o(n^h)$, where here and throughout the paper we
will use $\delta_H(p)$ to denote $p^{e(H)}(1-p)^{{h \choose
2}-e(H)}$. For a set of vertices $U \subseteq V$ we denote by
$H[U]$ the number of labeled copies of $H$ in $U$, and by $H^*[U]$
the number of induced labeled copies of $H$ in $U$. We also use
$e(U)$ to denote the number of edges inside a set of vertices $U$, and
$e(U,V)$ to denote the number of edges connecting $U$ to $V$. The following
is (part of) the main result of \cite{CGW}:

\begin{theo}[Chung, Graham, and Wilson \cite{CGW}]\label{theoCGW}
Fix any $0 < p < 1$. For any $n$-vertex graph $G$ the following
properties are equivalent:
\begin{itemize}
\item[${\cal P}_1$:] For any subset of vertices $U \subseteq
V(G)$ we have $~e(U)=\frac{1}{2}p|U|^2 +o(n^2)$.
\item[${\cal P}_2$:] For any subset of vertices $U \subseteq
V(G)$ of size $\frac12 n$ we have $~e(U)=\frac{1}{2}p|U|^2
+o(n^2)$.
\item[${\cal P}_3$:] Let $\lambda_i(G)$ denote the $i^{th}$
largest (in absolute value) eigenvalue of $G$. Then $e(G) =
\frac{1}{2}p n^2 +o(n^2)$, $\lambda_1(G)=pn +o(n)$ and
$\lambda_2(G)=o(n)$.
\item[${\cal P}_4(t)$:] For an even integer $t \geq 4$, let
$C_t$ denote the cycle of length $t$. Then $e(G) = \frac{1}{2}p
n^2 +o(n^2)$ and $C_t[G] = p^tn^t + o(n^t)$.
\item[${\cal P}_5$:] Fix an $\alpha \in (0,\frac12)$. For any
$U \subseteq V(G)$ of size $\alpha n$ we have $e(U, V \setminus
U)=p\alpha(1-\alpha)n^2 +o(n^2)$.
\end{itemize}
\end{theo}

For the rest of the paper let us use the notation $x=y \pm \epsilon$
as a shorthand for the two inequalities $y-\epsilon \leq x
\leq y + \epsilon$. As we have mentioned before, we use the $o(1)$ term to denote an arbitrary function tending to $0$ with
$n$. Hence, the meaning of the fact that, for example, ${\cal P}_2$ {\em
implies} ${\cal P}_1$ is that for any $f(n)=o(1)$ there is a $g(n)=o(1)$
such that if $G$ has the property that
all $U \subseteq V(G)$ of size $n/2$ satisfy
$e(U)=\frac{1}{2}p|U|^2 \pm  g(n)n^2$, then
$e(U)=\frac{1}{2}p|U|^2 \pm f(n)n^2$ for all $U \subseteq V(G)$.
Equivalently, this means for any $\delta>0$ there is an
$\epsilon=\epsilon(\delta)$ and $n_0=n_0(\delta)$ such that if $G$ is
a graph on $n \geq n_0$ vertices and it has the property that
all $U \subseteq V(G)$ of size $n/2$ satisfy
$e(U)=\frac{1}{2}p|U|^2 \pm \epsilon n^2$, then
$e(U)=\frac{1}{2}p|U|^2 \pm \delta n^2$ for all $U \subseteq V(G)$.
This will also be the meaning of other implications between
other graph properties later on in the paper.

%\footnote{An equivalent, more cumbersome, way to state Theorem
%\ref{theoCGW} would have been to replace all the $o(.)$ terms by
%$\epsilon_1,\ldots,\epsilon_5$ and say that there are some
%functions $f_{i,j}$ that relate $\epsilon_i$ and $\epsilon_j$.}

Note, that each of the items in Theorem \ref{theoCGW} is a property
we would expect $G(n,p)$ to satisfy with high probability. We will
thus say that $G$ is {\em $p$-quasi-random} if it satisfies property
${\cal P}_1$, that is if for some small $\delta$ all $U \subseteq
V(G)$ satisfy $e(U)=\frac{1}{2}p|U|^2 \pm \delta n^2$. If one wishes
to be more precise then one can in fact say that such a graph is
$(p,\delta)$-quasi-random. We will sometimes omit the $p$ and just
say that a graph is quasi-random. In the rest of the paper the
meaning of a statement ``If $G$ satisfies ${\cal P}_2$ then $G$ is
quasi-random'' is that ${\cal P}_2$ implies ${\cal P}_1$ in the
sense of Theorem \ref{theoCGW} discussed in the previous paragraph.
We will also say that a graph property ${\cal P}$ is {\em
quasi-random} if any graph that satisfies ${\cal P}$ must be
quasi-random. So the meaning of the statement ``${\cal P}_2$ is
quasi-random'' is that ${\cal P}_2$ implies ${\cal P}_1$. Therefore,
all the properties in Theorem \ref{theoCGW} are quasi-random.

Given Theorem \ref{theoCGW} one may think that any property
that holds with high probability in $G(n,p)$ is quasi-random. That
however, is far from true. For example, it is easy to see that
having the ``correct'' vertex degrees is not a quasi-random
property (consider $K_{n/2,n/2}$). Note also that in ${\cal P}_5$
we require $\alpha < \frac12$, because when $\alpha=\frac12$ the
property is not quasi-random (see \cite{CG3} and \cite{SSnind}). A
more relevant family of non quasi-random properties are those
requiring the graph to have the ``correct'' number of copies of a
fixed graph $H$. Note that ${\cal P}_4(t)$ guarantees that for any
even $t$, if a graph has the ``correct'' number of edges and the
``correct'' number of copies of $C_t$ then it is quasi-random. As
observed in \cite{CGW} this is {\em not} true for all graphs, in
fact this is not true for any non-bipartite $H$.

\subsection{Quasi-randomness and the distribution of copies of a single graph}

As throughout the paper we work with labeled copies and labeled
induced copies of $H$, we henceforth just call them copies and
induced copies. To understand the context of our main result, which
deals with induced copies of a fixed graph $H$, it is instructive
to review what is known about the effect of the distribution of a
fixed graph $H$ on quasi-randomness. By Theorem \ref{theoCGW} we
know that for some graphs $H$ the property of having the ``correct''
number of copies of $H$ in $G$, along with the right number of
edges, is enough to guarantee that $G$ is quasi-random.
Furthermore, this is not true for all graphs $H$. However, the
intuition is that something along these lines should be true for
any $H$, i.e. that for any $H$, if the copies of $H$ in a graph
$G$ have the ``properties'' we would expect them to have in
$G(n,p)$, then $G$ should be $p$-quasi-random. Simonovits and
S\'os \cite{SSnind} observed that the counter examples showing
that, for some graphs $H$, having just the ``correct'' number of copies
of $H$ (and the ``correct'' number of edges) is not enough to
guarantee quasi-randomness, all have the property that some of the
induced subgraphs of these counter examples have significantly
more/less copies of $H$ than we would expect to find in $G(n,p)$.
For example, in order to show that having the ``correct'' number of
edges and triangles as in $G(n,1/2)$ does not guarantee that $G$
is $\frac12$-quasi-random, one can take a complete graph on
$\alpha n$ vertices and a complete bipartite graph on
$(1-\alpha)n$ vertices, for an appropriate $\alpha$.

The main insight of Simonovits and S\'os \cite{SSnind} was that
quasi-randomness is a {\em hereditary} property, in the sense that
we expect a sub-structure of a random-like object to be
random-like as well. Thus, perhaps it will suffice to require that
the subgraphs of $G$ should also have the ``correct'' number of copies
of $H$. To state the main result of \cite{SSnind} let us introduce
the following variant of property ${\cal P}_1$ of Theorem
\ref{theoCGW}.

\begin{definition}[$H\mbox{[}U_1,\ldots,U_h\mbox{]}$]\label{DefCountCopies}
For a graph $H$ on $h$ vertices, and pairwise disjoint vertex sets $U_1,\ldots,U_h$,
we denote by $H[U_1,\ldots,U_h]$ the number of $h$-tuples $v_1 \in
U_{1},\ldots,v_h \in U_{h}$ that span a labeled copy of $H$.
\end{definition}

\begin{definition}[${\cal P}_H$] For a fixed graph $H$ on $h$
vertices, we say that a graph $G$ satisfies ${\cal P}_H$ if all
pairwise disjoint $h$-tuples $U_1,\ldots,U_h \subseteq V(G)$ of
equal (arbitrary) size $m$ satisfy
$$
H[U_1,\ldots,U_h]=p^{e(H)}h!m^h+o(n^h)\;.
$$
\end{definition}

Note that the above restriction is that the value of
$H[U_1,\ldots,U_h]$ should be close to what it should be in $G(n,p)$
for all $h$-tuples of equal-size. Observe also that the above
condition does not impose any restriction on the number of edges of
$G$, while in property ${\cal P}_1$ there is. Note also, that the
error in the above definition involves $n$ rather than
$m=|U_1|=\cdots=|U_h|$ so when $m=o(n)$ the condition vacuously
holds. As opposed to ${\cal P}_4$, which is not quasi-random for all
graphs, Simonovits and S\'os \cite{SSnind} showed that ${\cal P}_H$
is quasi-random for any graph $H$.

\begin{theo}[Simonovits and S\'os \cite{SSnind}]\label{theoSS}
The following holds for any graph $H$: if a graph $G$ satisfies
${\cal P}_H$ then it is $p$-quasi-random.
\end{theo}

Observe that ${\cal P}_H$ requires, via
Definition \ref{DefCountCopies}, all $h$-tuples of vertex sets to
have the ``correct'' number of copies of $H$ with one vertex in each
set. A more ``natural'' requirement, that was actually used in
\cite{SSnind}, is that all subsets of vertices $U \subseteq V(G)$
should contain the ``correct'' number of copies of $H$, that is, that
$H[U] \approx p^{e(H)}|U|^h$ for all $U \subseteq V(G)$. However,
it is not difficult to show that these two conditions are in fact
equivalent (see \cite{S}). We choose to work with Definition
\ref{DefCountCopies} as it will fit better with the discussion in
the next subsection.

\subsection{The main result}\label{subsecmain}

So we know from Theorem \ref{theoCGW} that when we consider the
number of subgraphs of $H$ in $G$, then some $H$ but not all, are
such that having the ``correct'' number of copies of $H$ in a graph
$G$ (and number of edges) is enough to guarantee that $G$ is
quasi-random. From Theorem \ref{theoSS} we know that for all $H$,
having the ``correct'' number of copies of $H$ in all the subgraphs of
$G$ is enough to guarantee that $G$ is quasi-random. A natural
question is what can we learn from the distribution of {\em
induced} copies of a graph $H$? As we shall see, the situation is much
more involved.

Recall that for a fixed graph $H$ on $h$ vertices and a fixed $0<
p < 1$, we define $\delta_H(p)=p^{e(H)}(1-p)^{{h \choose
2}-e(H)}$. Let us denote by $\overline{p}_H$ the second
\footnote{It is not difficult to see that for non-negative
integers $k,\ell$ the equation $x^k(1-x)^\ell=q$ has at most two
solutions in $(0,1)$.} solution (other than $p$) of the equation
$\delta_H(p)=x^{e(H)}(1-x)^{{h \choose 2}-e(H)}$. We call
$\overline{p}_H$ the {\em conjugate} of $p$ with respect to $H$.
We will sometimes just write $\overline{p}$ instead of
$\overline{p}_H$ when $H$ is fixed. Note that the expected number
of induced copies of $H$ in a set of vertices $U$ is roughly
$\delta_H(p)|U|^h$. But, as it may \footnote{The only case where
$p=\overline{p}_H$ is when $p=e(H)/{h \choose 2}$ or when $e(H) \in \{0,{h \choose 2}\}$.} be the case
that $p \neq \overline{p}_H$ we see that for any $H$ and any $p$,
the distribution of induced copies of $H$ in both $G(n,p)$ and
$G(n,\overline{p}_H)$ behaves {\em precisely} the same. Therefore,
the best we can hope to deduce from the fact that the distribution
of induced copies of $H$ in $G$ is close to that of $G(n,p)$ is
that $G$ is either $p$-quasi-random or
$\overline{p}_H$-quasi-random.

Let us denote by $H^*[U_1,\ldots,U_h]$ the natural generalization
of $H[U_1,\ldots,U_h]$ (defined in Definition
\ref{DefCountCopies}) with respect to induced subgraphs, that is,
$H^*[U_1,\ldots,U_h]$ is the number of $h$-tuples of vertices $v_1
\in U_{1},\ldots,v_h \in U_{h}$ with the property that
$v_1,\ldots,v_h$ span a labeled induced copy of $H$. Note that for
an $h$ tuple of vertex sets $U_1,\ldots,U_h$ in $G(n,p)$ each of
size $m$, the expected value of $H^*[U_1,\ldots,U_h]$ is
$\delta_H(p)h! m^h$.

So given the above discussion and Theorem \ref{theoSS}, it seems
reasonable to conjecture that, if a graph $G$ has the ``correct''
distribution of induced copies of $H$, then $G$ is either
$p$-quasi-random or $\overline{p}_H$-quasi-random. When we say
``correct'' distribution we mean that all pairwise disjoint $h$-tuples
$U_1,\ldots,U_h \subseteq V(G)$ of the same size $m$ satisfy
$H^*[U_1,\ldots,U_h] = \delta_H(p)h! m^h \pm o(m^h)$. However, it was
observed in \cite{SSnind,SSind} that this is not the case. For
example, one can take vertex set $V_1,V_2$ of sizes $\alpha n,
(1-\alpha) n$ and put $G(\alpha n,p_1)$ on $V_1$,
$G((1-\alpha_1)n,p_1)$ on $V_2$ and connect $V_1$ and $V_2$ with
probability $p_2 \neq p_1$. Then for appropriate constants, we get
a graph with the ``correct'' distribution of the 3-vertex path, yet
this graph is not $p$-quasi-random for any $p$.

However, as before, the intuition is that having the ``correct''
distribution of induced copies of $H$ {\em should} guarantee that
$G$ is quasi-random. Our main result in this paper is that indeed
it does, one just needs to refine the notion of ``correct
distribution''. As we have mentioned before, if $U_1,\ldots,U_h$
is an $h$-tuple of vertices in $G(n,p)$ of the same size $m$, then
we would expect to have $H^*[U_1,\ldots,U_h] = \delta_H(p)h!
m^h \pm o(m^h)$. However, this is because we would actually expect a slightly stronger
condition to hold. Before
stating this condition, let us introduce the following
``permuted'' version of the quantity $H^*[U_1,\ldots,U_h]$.

\begin{definition}[$H^*_{\sigma}\mbox{[}U_1,\ldots,U_h\mbox{]}$]
Let $H$ be a graph on $h$ vertices, let $U_1,\ldots,U_h$ be an
$h$-tuple of pairwise disjoint vertex sets, and let $\sigma \in S_h$ be a permutation
$[h] \rightarrow [h]$. Then we denote by
$H^*_{\sigma}[U_1,\ldots,U_h]$ the number of $h$-tuples of
vertices $v_1 \in U_{\sigma(1)},\ldots,v_h \in U_{\sigma(h)}$ with
the property that $v_i \in U_{\sigma(i)}$ is connected to $v_j \in
U_{\sigma(j)}$ if and only if $(i,j) \in H$.
\end{definition}

Getting back to our discussion, observe that the reason we expect
to have $H^*[U_1,\ldots,U_h] = \delta_H(p)h! m^h \pm o(m^h)$ is simply
because we expect to have $H^*_{\sigma}[U_1,\ldots,U_h] =
\delta_H(p)m^h \pm o(m^h)$ for all $h!$ permutations in $S_h$. It is now
natural to define the following property:

\begin{definition}[${\cal P}^*_H$] For a fixed graph $H$ on $h$
vertices, we say that a graph $G$ satisfies ${\cal P}^*_H$ if for
all pairwise disjoint $h$-tuples $U_1,\ldots,U_h \subseteq V(G)$
of equal (arbitrary) size $m$, and for every $\sigma \in S_h$
$$
H^*_{\sigma}[U_1,\ldots,U_h] = \delta_H(p)m^h + o(n^h)\;.
$$
\end{definition}

Our main result is that property ${\cal P}^*_H$ guarantees that a
graph is quasi-random.

\begin{theo}[Main result]\label{theomain} The following holds for
any graph $H$: if a graph satisfies ${\cal P}^*_H$ then it is
either $p$-quasi-random or $\overline{p}_H$-quasi-random.
\end{theo}

Our main result can be formulated as saying that for any $H$, if a
graph $G$ has the ``correct'' distribution of induced copies of $H$,
then $G$ is quasi-random. We remind the reader that one cannot
hope to strengthen Theorem \ref{theomain} by showing that $G$ must
be $p$-quasi-random, as $G(n,\overline{p})$ satisfies ${\cal
P}^*_H$ with probability 1. Observe, that our notion of ``correct
distribution'' (that is, the quantities $H^*_{\sigma}$) is stronger
than the notions that have been considered
before (that is, the quantities $H^*$), where the latter is known
to be too weak to guarantee quasi-randomness for arbitrary graphs $H$.
Let us mention here that Simonovits and S\'os conjectured in \cite{SSnind} that
the weaker quantities $H^*$ should be sufficient for guaranteeing quasi-randomness
for any $H$ on at least 4 vertices. This conjecture, however, is still wide open.

\subsection{Overview of the paper}

As we have discussed in the first subsection, the theory of
quasi-random graphs has many applications in theoretical computer
science, both in the case of sparse and dense graphs.
We think that the main interest of
our result is in the proof techniques and tools that are used in
the course of its proof. Besides several combinatorial arguments
and tools (such as the Regularity Lemma \cite{Sz}, Ramsey's
Theorem and R\"odl's ``nibble'' Theorem \cite{R}) the main
underlying idea of the proof is an algebraic one. Roughly
speaking, what we do is take all the information we know about the
graph $G$, namely the information on the distribution of induced
copies of $H$, and use it in order define a large system of
polynomial equations. The unknowns in this system of equations
represent (in some way) the distribution of edges of $G$. The crux
of the proof is to show that the unique solution of this system of
equations, is one that forces the edges of the graph to be nicely
distributed (in the sense of property ${\cal P}_1$ in Theorem
\ref{theoCGW}). The main theorem we need in order to obtain this
uniqueness is a result of Gottlieb \cite{G}, in algebraic
combinatorics, concerning the rank of set inclusion matrices (see
Theorem \ref{gottlieb}). This approach to showing that a graph is
quasi-random may be applicable for showing quasi-random properties
of other structures.

In Section \ref{SecMain} we prove Theorem \ref{theomain} by
applying several combinatorial tools as well as a key lemma (Lemma
\ref{epsilon}) that is proved in Section \ref{SecKey}. The proof
of Lemma \ref{epsilon}, which is the most difficult step in the
proof of Theorem \ref{theomain}, contains most of the new ideas we
introduce in this paper.

\section{Proof of Main Result}\label{SecMain}

In this section we give the proof of Theorem \ref{theomain}, but
before getting to the actual proof we will need some preparation. We
first discuss Lemma \ref{epsilon}, which is the main technical lemma
we need for the proof of Theorem \ref{theomain}, and whose proof
appears in the next section. We then discuss some simple notions
related to the Regularity Lemma, and then turn to the proof of
Theorem \ref{theomain}. Throughout this section, let us fix a real
$0 < p < 1$ and a graph $H$ on $h$ vertices. Recall that we set
$\delta_H(p)=p^{e(H)}(1-p)^{{h \choose 2}-e(H)}$ and that we denote
by $\overline{p}$, the conjugate of $p$, the second solution in
$(0,1)$ of the equation $\delta_H(p)=x^{e(H)}(1-x)^{{h \choose
2}-e(H)}$.

\subsection{The Key Lemma}

In what follows we will work with weighted complete graphs $W$ on
$r$ vertices. We will think of the vertices of $W$ as the integers
$[r]$. In that case each pair of vertices $1 \leq i < j \leq r$
will have a weight $0 \leq w(i,j) \leq 1$. Let us identify the $h$
vertices of $H$ with the integers $[h]$. Given an injective
mapping $\phi:[h] \rightarrow [r]$, which we think of as a mapping
from the vertices of $H$ to the vertices of $W$, we will set
$$
W(\phi) = \prod_{(i,j) \in E(H)}w(\phi(i),\phi(j))\prod_{(i,j)
\not \in E(H)}(1-w(\phi(i),\phi(j)))\;.
$$
Another notation that will simplify the presentation is a variant
of the $H^*_{\sigma}[U_1,\ldots,U_h]$ notation that was defined in
Section \ref{SecIntro}. Suppose we have $r$ pairwise disjoint
vertex sets $U_1,\ldots,U_r$ and an injective mapping $\phi:[h]
\rightarrow [r]$. Then we denote by $H^*_{\phi}[U_1,\ldots,U_r]$
the number of $h$-tuples of vertices $v_1 \in
U_{\phi(1)},\ldots,v_h \in U_{\phi(h)}$ with the property that
$v_i \in U_{\phi(i)}$ is connected to $v_j \in U_{\phi(j)}$ if and
only if $(i,j)$ is an edge of $H$.

Suppose we construct an $r$-partite graph on vertex sets
$U_1,\ldots,U_r$, each of size $m$, by connecting every vertex in
$U_i$ with any vertex in $U_j$ independently with probability
$w(i,j)$. Then, observe that for any $\phi: [h] \rightarrow [r]$,
we would expect $H^*_{\phi}[U_1,\ldots,U_r]$ to be close to
$W(\phi)m^h$. Continuing this example, suppose that all $(i,j)$
satisfy $w(i,j)=p$. Then we would expect all $\phi$ to satisfy
$H^*_{\phi}[U_1,\ldots,U_h] = \delta_H(p)m^h$. Observe however,
that we would also expect the same to hold if we were to replace
$p$ by $\overline{p}$.

The following lemma shows that the converse is also true in the
following sense: if we know that for {\em any} injective mapping
$\phi$ we have the ``correct'' fraction of induced copies of $H$ as we
would expect to find if we had $w(i,j)=p$ for all $(i,j)$, then
either\footnote{Remember that we cannot expect to be able to show
that all densities are $p$ as the number of induced copies of $H$
behaves the same with respect to $p$ and $\overline{p}$.} almost
all $(i,j)$ satisfy $w(i,j)=p$ or almost all satisfy
$w(i,j)=\overline{p}$. Note that for convenience the lemma is
stated with respect to quantities in $(0,1)$, rather than with
respect to the number of edges or number of copies \footnote{So
$w(i,j)$ should be understood as the density between the pair
$(U_i,U_j)$ and $W(\phi)$ is $H^*_{\phi}[U_1,\ldots,U_r]/m^h$.}.
In what follows, we will always assume wlog that if $p \neq
\overline{p}$ then $\epsilon< |p-\overline{p}|/2$. This will
guarantee that the intervals $\overline{p} \pm \epsilon$ and $p \pm \epsilon$
are disjoint.

\begin{lemma}[The Key Lemma]
\label{epsilon} For every $h$ there exists an $N_{\ref{epsilon}}=
N_{\ref{epsilon}}(h)$ so that for any $r \geq N_{\ref{epsilon}}$
and $\epsilon > 0$ there exists
$\delta_{\ref{epsilon}}=\delta_{\ref{epsilon}}(\epsilon,h,r)
> 0$ with the following properties: suppose $W$ is a weighted graph
on $r$ vertices, such that for all $\phi: [h] \rightarrow [r]$ we
have $W(\phi) = \delta_H(p) \pm \delta_{\ref{epsilon}}$. Then any
pair $(i,j)$ satisfies either $w(i,j) = p \pm \epsilon$ or $w(i,j)
= \overline{p} \pm \epsilon$. Furthermore, either at most $r-1$ of
the pairs $(i,j)$ satisfy $w(i,j) = p \pm \epsilon$ or at most
$r-1$ of the pairs $(i,j)$ satisfy $w(i,j) = \overline{p} \pm
\epsilon$.
\end{lemma}

The proof of Lemma \ref{epsilon}, which is the main lemma we need
for the proof of Theorem \ref{theomain}, appears in Section
\ref{SecKey}. It is interesting to note that, as we show in Section
\ref{SecKey}, one cannot strengthen the above lemma by showing
that either all densities are close to $p$ or they are all close
to $\overline{p}$.

\subsection{The Regularity Lemma}

We now give a brief overview of the Regularity Lemma of
Szemer\'edi, which turns out to be strongly related to
quasi-random graphs. For a pair of nonempty vertex sets $(A,B)$ we
denote by $d(A,B)$ the edge density between $A$ and $B$, that is
$d(A,B)=|E(A,B)|/|A||B|$. A pair of vertex sets $(A,B)$ is said to
be {\em $\gamma$-regular}, if for any two subsets $A' \subseteq A$
and $B' \subseteq B$, satisfying $|A'| \geq \gamma|A|$ and $|B'|
\geq \gamma|B|$, the inequality $|d(A',B')-d(A,B)| \leq \gamma$
holds. A partition of the vertex set of a graph is called an {\em
equipartition} if all the sets of the partition are of the same
size (up to $1$). We call the number of partition classes of an
equipartition the {\em order} of the equipartition. Finally, an
equipartition ${\cal V} = \{V_1,\ldots,V_k \}$ of the vertex set
of a graph is called {\em $\gamma$-regular} if all but at most
$\gamma{k \choose 2}$ of the pairs $(V_i,V_{j})$ are
$\gamma$-regular. The celebrated Regularity Lemma of Szemer\'edi
can be formulated as follows:

\begin{lemma}[\cite{Sz}]\label{RegLemma}
For every $t$ and $\gamma > 0$ there exists $T =
T_{\ref{RegLemma}}(\gamma,t)$, such that any graph of size at
least $t$ has a $\gamma$-regular equipartition of order $k$, where
$t \leq k \leq T$.
\end{lemma}

The following lemma of Simonovits and S\'os \cite{SS} shows that
the property of having a regular partition where most of the pairs
are connected by regular pairs of density close to $p$ implies
that the graph is $p$-quasi-random. For completeness we include a
short self-contained proof of this lemma at the end of this
section. We remind the reader that we use the notation $x=y \pm \epsilon$ to denote
the fact that $y - \epsilon \leq x \leq y +\epsilon$.

\begin{lemma}[Simonovits and S\'os \cite{SS}]\label{SSlemma}
For every $\zeta >0$ there is an
$\epsilon=\epsilon_{\ref{SSlemma}}(\zeta)$ and
$t=t_{\ref{SSlemma}}(\zeta)$ with the following property: suppose
an $n$ vertex graph $G$ has an $\epsilon$-regular partition of
order $k \ge t$ where all but $\epsilon {k \choose 2}$ of the
pairs are $\epsilon$-regular with density $p \pm \epsilon$. Then
every set of vertices $U \subseteq G$ spans $\frac12p|U|^2 \pm
\zeta n^2$ edges.
\end{lemma}

Another tool we will need for the proof of Theorem \ref{theomain}
is Lemma \ref{cntHWsig} below.
\ignore{Note that in this lemma it is enough to consider only the
case $r=h$. We have stated with respect to $r$ sets of vertices in
order to simplify its application in the proof of Theorem
\ref{theomain}.} This lemma is equivalent to saying that if we
have $r$ pairwise disjoint vertex sets $V_1,\ldots,V_r$ that are
all regular enough, then for any injective mapping $\phi:[h]
\rightarrow [r]$ we have that $H^*_{\phi}[V_1,\ldots,V_r]$ is
close to what it should be. Such a lemma is well known, and has
been proven and used in many papers. See, e.g., Lemma 4.2 in
\cite{F} for one such proof. We thus omit the proof of Lemma
\ref{cntHWsig}.

\begin{lemma}\label{cntHWsig} For any $\delta>0$ and $h$, there
exists a $\gamma=\gamma_{\ref{cntHWsig}}(\delta,h)>0$ such that
the following holds: Let $W$ be a weighted complete graph on $r$
vertices, and suppose $V_1,\ldots,V_r$ are pairwise disjoint sets
of vertices of size $m$ each, that all pairs $(V_i,V_j)$ are
$\gamma$-regular and that all pairs satisfy $d(V_i,V_j)=w(i,j)$.
Then, for any injective mapping $\phi:[h] \rightarrow [r]$, we
have
\begin{equation}\label{wisb}
H_{\phi}^*[V_1,\ldots,V_r] = (W(\phi) \pm \delta) m^h.
\end{equation}
\end{lemma}

\subsection{Proof of Theorem \ref{theomain}}

For the proof of Theorem \ref{theomain} we will also need the
following two lemmas, whose proofs are deferred to the end of this
section.

\begin{lemma}\label{estdense}
For every $\epsilon$ there is an
$r_{\ref{estdense}}=r_{\ref{estdense}}(\epsilon)$ such that for
every $r \geq r_{\ref{estdense}}$ there is
$N_{\ref{estdense}}=N_{\ref{estdense}}(r)$ and
$\gamma_{\ref{estdense}}=\gamma_{\ref{estdense}}(r)$ with the
following property. Assume $k \geq N_{\ref{estdense}}$ and that
$K$ is a $k$ vertex graph with at least
$(1-\gamma_{\ref{estdense}}){k \choose 2}$ edges. Suppose the
edges of $K$ are colored red/blue so that at least $\epsilon {k
\choose 2}$ are blue and at least $\epsilon {k \choose 2}$ are
red. Then $K$ has $r$ vertices that span a complete graph $K_r$
with at least $r$ blue edges and at least $r$ red edges.
\end{lemma}

\begin{lemma}\label{coveredges} For any $\gamma$ and $r$, there is an
$N_{\ref{coveredges}}= N_{\ref{coveredges}}(\gamma,r)$ such that
the following holds for any $k \geq N_{\ref{coveredges}}$. If $K$
is a graph on $k$ vertices with at least $(1-\gamma){k \choose 2}$
edges, then $K$ has at least $(1-\gamma r^2){k \choose 2}$ edges
that belong to a copy of $K_r$.
\end{lemma}

\paragraph{Proof of Theorem \ref{theomain}:} We will say that a
$\gamma$-regular equipartition of order $k$ is
$\gamma$-super-regular\footnote{In some papers the term $\gamma$-super-regular is used
for other notions of regular partitions of graphs, but these should not be confused with
the notion we define and use here.} if all but $\gamma {k \choose 2}$ of the
pairs are $\gamma$-regular with density $p \pm \gamma$ or all but
$\gamma {k \choose 2}$ of the pairs are $\gamma$-regular with
density $\overline{p} \pm \gamma$. We need to show that any graph
$G$ that satisfies ${\cal P}^*_{H}$ must be either
$p$-quasi-random or $\overline{p}$-quasi-random\footnote{Recall
that the meaning of that is that either every set $U \subseteq
V(G)$ satisfies $e(U)=\frac12p |U^2| \pm \zeta n^2$ or that every
such set satisfies $e(U)=\frac12\overline{p} |U^2| \pm \zeta n^2$
for some small $\zeta
>0$}. Fix any $\zeta
>0$ and recall that by Lemma \ref{SSlemma} we know that in order
to show that a graph $G$ has the property that every set $U
\subseteq V(G)$ satisfies $e(U)=\frac12p |U^2| \pm \zeta n^2$ or
that every such set satisfies $e(U)=\frac12\overline{p} |U^2| \pm
\zeta n^2$, it is enough to show that $G$ has an
$\epsilon$-super-regular partition of order at least $t$, where
\begin{equation}\label{chooset}
t=t_{\ref{SSlemma}}(\zeta)\;,
\end{equation}
and
\begin{equation}\label{chooseeps}
\epsilon=\epsilon_{\ref{SSlemma}}(\zeta)\;.
\end{equation}
Let us define the following constants\footnote{We note that we
need $N$, which is defined in (\ref{chooseN}), in order to allow
us to apply the various lemmas we stated above, that all work for
large enough graphs.}
\begin{equation}\label{chooser}
r=\max(r_{\ref{estdense}}(\epsilon/2),~N_{\ref{epsilon}}(h))
\end{equation}
\begin{equation}\label{choosegamma}
\gamma=\min(\epsilon,~\gamma_{\ref{estdense}}(r)~,\gamma_{\ref{cntHWsig}}(\delta_{\ref{epsilon}}(\epsilon,h,r)/2,h))
\end{equation}
\begin{equation}\label{chooseN}
N=\max(t,~N_{\ref{coveredges}}(\gamma,r),~N_{\ref{estdense}}(r),~N_{\ref{epsilon}}(h))
\end{equation}
\begin{equation}\label{chooseT}
T=T_{\ref{RegLemma}}(\gamma/r^2,~ N)
\end{equation}
\begin{equation}\label{choosedelta}
\delta=\delta_{\ref{epsilon}}(\epsilon,h,r)/2T^h.
\end{equation}
To complete the proof that ${\cal P}^*_{H}$ implies that a graph
is either $p$-quasi-random or $\overline{p}$-quasi-random, we show
(via Lemma \ref{SSlemma}) that for any $\zeta >0$ there is an
$N(\zeta)$ and $\delta(\zeta)$ such that the following holds: if
$G$ is a graph on at least $N(\zeta)$ vertices and for every
$h$-tuple of vertex sets $U_1,\ldots,U_h \subseteq V(G)$ of
(arbitrary) size $m$ each, and for every permutation $\sigma: [h] \rightarrow
[h]$ we have
\begin{equation}\label{condcopies}
H_{\sigma}^*[U_1,\ldots,U_h]=\delta_H(p) m^h \pm \delta(\zeta)
n^h\;,
\end{equation}
then $G$  has an $\epsilon$-super-regular partition of order $k$,
where $k \geq t$. We will show that one can take $N(\zeta)$ to be
the integer $N$ defined in (\ref{chooseN}) and that
$\delta(\zeta)$ can be taken as the value defined in
(\ref{choosedelta}).

So let $G$ be a graph of size at least $N$, and apply Lemma
\ref{RegLemma} (the regularity lemma) on $G$ with $\gamma/r^2$ and
$N$ that were defined in (\ref{choosegamma}) and (\ref{chooseN}).
Lemma \ref{RegLemma} guarantees that $G$ has a $\gamma$-regular
partition ${\cal V}=\{V_1,\ldots,V_k\}$ where $t \leq N \leq k
\leq T$ and $T$ is given in (\ref{chooseT}). We now need to show
that all but $\epsilon {k \choose 2}$ of the pairs $(V_i,V_j)$ are
$\epsilon$-regular and satisfy $d(V_i,V_j)=p \pm \epsilon$ or
$\epsilon$-regular and satisfy $d(V_i,V_j)=\overline{p} \pm
\epsilon$. Let us define $W$ to be a weighted graph on $k$
vertices, where if $(V_i,V_j)$ is $\epsilon$-regular then $(i,j)$
are connected with an edge of weight $w(i,j)=d(V_i,V_j)$, and if
$(V_i,V_j)$ is not $\epsilon$-regular then $(i,j)$ are not
connected. So our goal is to show that either all but $\epsilon {k
\choose 2}$ of the pairs of vertices of $W$ are connected by an
edge with weight $p \pm \epsilon$ or that all but $\epsilon {k
\choose 2}$ of the pairs of vertices of $W$ are connected by an
edge with weight $\overline{p} \pm \epsilon$.

\begin{claim}\label{noKr} Any copy of $K_r$ in $W$ satisfies the
following:
\begin{enumerate}
\item Any edge has either weight $p \pm \epsilon$ or weight
$\overline{p} \pm \epsilon$.
\item If $p \neq \overline{p}$ then either at most $r-1$ of them
have weight $p \pm \epsilon$ or at most $r-1$ of them have weight
$\overline{p} \pm \epsilon$.
\end{enumerate}
\end{claim}

\paragraph{Proof:} Consider any copy of $K_r$ in $W$ and suppose
wlog that the vertices of this copy are $1,\ldots,r$. Recall that
$k \leq T$ and that $G$ satisfies (\ref{condcopies}) with the
$\delta$ that was chosen in (\ref{choosedelta}). We thus infer
that for any injective mapping $\phi:[h] \rightarrow [r]$
\begin{equation}\label{numcopies1}
H^*_{\phi}[V_{1},\ldots,V_{r}]=\delta_H(p)\left(\frac{n}{k}\right)^h
\pm \delta n^h=(\delta_H(p) \pm
\frac12\delta_{\ref{epsilon}}(\epsilon,h,r))\left(\frac{n}{k}\right)^h\;.
\end{equation}
In addition, as we are referring to $r$ vertices that form a copy
of $K_r$ in $W$, we know that $V_{1},\ldots,V_{r}$ are all
pairwise $\gamma$-regular. Thus the choice of $\gamma$ in
(\ref{choosegamma}) guarantees via Lemma \ref{cntHWsig} that for
any injective mapping $\phi:[h] \rightarrow [r]$ we have
\begin{equation}\label{numcopies2}
H^*_{\phi}[V_{1},\ldots,V_{r}]=(W(\phi) \pm
\frac12\delta_{\ref{epsilon}}(\epsilon,h,r))\left(\frac{n}{k}\right)^h\;.
\end{equation}
Combining (\ref{numcopies1}) and (\ref{numcopies2}) we infer that
for any $\phi:[h] \rightarrow [r]$ we have $W(\phi)=\delta_H(p)
\pm \delta_{\ref{epsilon}}(\epsilon,h,r)$. Hence, the two
assertions of the claim follow from Lemma \ref{epsilon}. $\qed$

\bigskip

We are now going to use Lemma \ref{epsilon} in order to color some
of the edges of $W$. Consider any copy of $K_r$ in $W$. We know
from the first assertion of Claim \ref{noKr} that all the edge
weights in the copy of $K_r$ are either $p \pm \epsilon$ or
$\overline{p} \pm \epsilon$. If $p=\overline{p}$ then we color all
the edges of this $K_r$ with the color red. So assume that $p \neq
\overline{p}$ and recall that we assume wlog that in this case
$\epsilon < |p-\overline{p}|/2$, which makes it possible to color
the edges whose weight is $p \pm \epsilon$ blue, and the edges
whose weight is $\overline{p} \pm \epsilon$ red (in a well defined
way). We now apply this coloring scheme to any copy of $K_r$ in
$W$. We claim that we have thus colored at least $(1-\gamma){k
\choose 2}$ of the edges of $W$. Indeed, as we applied the
regularity lemma with $\gamma/r^2$ we know that $W$ has at least
$(1-\gamma/r^2){k \choose 2}$ edges. As $k \geq
N_{\ref{coveredges}}(\gamma,r)$ we infer from Lemma
\ref{coveredges} that at least $(1-\gamma){k \choose 2}$ of the
edges of $W$ belong to a copy of $K_r$ thus they are colored in
the above process. Let us now remove from $W$ all the uncolored
edges and call the new graph $W'$. Thus $W'$ has at least
$(1-\gamma){k \choose 2}$ edges and they are all colored either
red or blue.

We now claim that either $W'$ has at most $\frac{\epsilon}{2} {k
\choose 2}$ red edges, or at most $\frac{\epsilon}{2}{k \choose
2}$ blue edges. Indeed, if $W'$ had at least $\epsilon {k \choose
2}$ red edges {\em and} at least $\epsilon {k \choose 2}$ blue
edges, then our choice of $r$ and $\gamma$ in (\ref{chooser}) and
(\ref{choosegamma}), the fact that $W'$ has at least $(1-\gamma){k
\choose 2}$ edges, and that $k \geq N_{\ref{estdense}}(r)$, would
allow us to apply Lemma \ref{estdense} on $W'$ and infer that it
has a copy of $K_r$ with at least $r$ blue edges and at least $r$
red edges, contradicting Claim \ref{noKr} (recall that $W'$ is a
subgraph of $W$).

We thus conclude that $W$ has at least $(1-\gamma){k \choose 2}
\geq (1-\frac{\epsilon}{2}){k \choose 2}$ edges, and that even if
$p \neq \overline{p}$ either all but $\frac{\epsilon}{2}{k \choose
2}$ of them are red or all but at most $\frac{\epsilon}{2}{k
\choose 2}$ of them are blue. By the definition of $W$, this means
that in the equipartition ${\cal V}$ either all but $\epsilon {k
\choose 2}$ of the pairs are $\epsilon$-regular with density $p
\pm \epsilon$, or all but at most $\epsilon {k \choose 2}$ of them
are $\epsilon$-regular with density $\overline{p} \pm \epsilon$,
which completes the proof. $\qed$

\subsection{Proofs of additional lemmas}

We end this section with the proofs of Lemmas \ref{SSlemma},
\ref{estdense} and \ref{coveredges}.

\paragraph{Proof of Lemma \ref{SSlemma}:} We claim that one can take
$\epsilon=\epsilon_{\ref{SSlemma}}(\zeta)=\frac18\zeta$ and
$t=t_{\ref{SSlemma}}(\zeta)=8/\zeta$. Indeed, suppose $G$ has an
$\epsilon$-super-regular partition $\{V_1,\ldots,V_k\}$ of order $k
\geq t$, that is, a partition in which all but $\epsilon {k \choose
2}$ of the pairs $(V_i,V_j)$ are $\epsilon$-regular with density
$d(V_i,V_j)=p \pm \epsilon$. Let us count the number of edges of $G$
that do not connect a pair $(V_i,V_j)$ which is $\epsilon$-regular
with density $d(V_i,V_j)=p \pm \epsilon$. As $k \geq t \geq 8/\zeta$
we know that the number of pairs of vertices that both belong to the
same set $V_i$ is at most $k |n/k|^2 \leq \frac18 \zeta n^2$. As all
but $\epsilon {k \choose 2}$ of the pairs $(V_i,V_j)$ are
$\epsilon$-regular with density $d(V_i,V_j)=p \pm \epsilon$, we also
know that the number of pairs connecting pairs $(V_i,V_j)$, which do
not satisfy these two conditions, is bounded by $\epsilon {k \choose
2}(n/k)^2 \leq \frac18 \zeta n^2$.

Consider now a set of vertices $U$, and define $U_i=U \cap V_i$. The
number of vertices of $U$ that belong to a set $U_i$ whose size is
smaller than $\epsilon |V_i|$ is bounded by $\epsilon n$. Therefore
the number of pairs of vertices of $U$ such that one of them belongs
to a set $U_i$ of size smaller than $\epsilon |V_i|$ is bounded by
$\epsilon n^2 \leq \frac18 \zeta n^2$. Combining the above three
facts we conclude that all but $\frac12\zeta n^2$ of the pairs of
vertices of $u_i,u_j \in U$ are such that: (1) $u_i \in U_i$, $u_j
\in U_j$ and $i \neq j$; (2) $U_i \geq \epsilon |V_i|$ and $U_j \geq
\epsilon |V_j|$; (3) $(V_i,V_j)$ is $\epsilon$-regular with density
$p \pm \epsilon$. Therefore, by the definition of a regular pair we
get that the density of $U$ in all but $\frac12\zeta n^2$ of its
pairs is $p \pm 2 \epsilon = p \pm \frac12\zeta$, and therefore
$e(U)=\frac12p|U|^2 \pm \zeta n^2$. $\qed$

\paragraph{Proof of Lemma \ref{estdense}:} Suppose we randomly pick
$r$ vertices $v_1, \ldots, v_r$ from $K$ {\em with repetitions}
where $r=\Omega(1/\epsilon^3)$. Clearly, if $k \geq 10r^2$ then by a
Birthday-Paradox argument we infer that with probability at least
$3/4$ all the vertices $v_1, \ldots, v_r$ are distinct. Suppose wlog
that $r$ is even and let us partition the set of unordered pairs
$(v_i,v_j)$ into $r-1$ perfect matchings $M_1,\ldots,M_{r-1}$ on the
vertices $v_1,\ldots,v_r$. For every pair $(v_i,v_j)$ let $p_{i,j}$ be the
indicator random variable for the event that $v_i$ and $v_j$ are
connected in $K$ by a red edge. As we sample with repetitions then
for every matching $M_t$, the $r/2$ events $\{p_{i,j}: (v_i,v_j) \in M_t
\}$ are independent. Also, as $K$ has at least $\epsilon{k \choose
2}$ red edges, we have that $Pr[p_{i,j}=1] \geq \epsilon/2$ (we lose
a little due to the probability of having non distinct vertices). We
thus conclude that for any matching $M_t$, the expected number of
red edges spanned by its members is at
least $\epsilon r/4$ and by a Chernoff bound, the probability of
deviating from this expectation by more than $\epsilon r /8$ is
bounded by $2^{-\Theta(\epsilon^2 r)}<1/4r$. Clearly the same
analysis applies for the blue edges. We conclude by the union bound
that with probability at least 3/4 the $r$ vertices span at least
$\epsilon r^2/16$ red edges and at least $\epsilon r^2/16$ blue
edges. As $r=\Omega(1/\epsilon^3)$ we have $\epsilon r^2/16 \geq r$
therefore we have the required amount of red/blue edges. We conclude
that one can take $r_{\ref{estdense}}=\Omega(1/\epsilon^3)$ and
$N_{\ref{estdense}}(r)=10 r^2$.

Finally, to conclude that all the pairs $(v_i, v_j)$ are connected
we take $\gamma_{\ref{estdense}}=1/4r^2$. This way, the probability
that a pair of vertices are not connected is at most $1/4r^2$ and by
the union bound, with probability at least $3/4$ they are all
connected. So to recap, if we sample with repetition $r$ vertices,
then with probability at least $1/4$ they are all distinct, all
connected, and have at least $r$ red edges and at least $r$ blue
edges. So there must be at least one such set of $r$ vertices in
$K$. $\qed$

\paragraph{Proof of Lemma \ref{coveredges}:} Suppose $k$ is large
enough to guarantee by R\"odl's theorem \cite{R} that the complete
graph on $k$ vertices contains $(1-\gamma){k \choose 2}/{r \choose
2}$ edge disjoint copies of $K_r$. If we now consider the same
copies of $K_r$ in $K$ (more precisely, the vertex sets of these
copies) then the fact that $K$ has $(1-\gamma){k \choose 2}$ edges
implies that at most $\gamma {k \choose 2} = \gamma {r \choose 2}
\cdot {k \choose 2}/{r \choose 2}$ of these copies of $K_r$ have a
pair of vertices that are not connected. Thus, $K$ contains at least
$(1-\gamma r^2){k \choose 2}/{r \choose 2}$ edge disjoint copies of
$K_r$ implying that at least $(1-\gamma r^2){k \choose 2}$ edges of
$K$ belong to a copy of $K_r$. $\qed$

\section{Proof of the Key Lemma}\label{SecKey}

As in Section \ref{SecMain}, let us fix a real $0 < p < 1$ and a
fixed graph $H$ on $h$ vertices. Let also $\overline{p}$ be the
conjugate of $p$ with respect to $H$. We will again work with
weighted complete graphs $W$ on $r$ vertices, and will identify
the vertices of $W$ with $[r]$ and the vertices of $H$ with $[h]$.
Each pair of vertices $1 \leq i < j \leq r$ of $W$ has a weight $0
\leq w(i,j) \leq 1$ that is given by some weight function $w:E(W)
\rightarrow [0,1] $. We remind the reader of the notation
$W(\phi)$ that was introduced at  the beginning of Section
\ref{SecMain}.

Recall that Lemma \ref{epsilon} states that, if all the values
$W(\phi)$ are close to what they should be, then all the weights
$w(i,j)$ are close to what they should be. The following lemma is
an ``exact'' version of Lemma \ref{epsilon} where we assume that
the values $W(\phi)$ are exactly what they should be. The proof of
Lemma \ref{epsilon} will then follow from the lemma below using
standard continuity arguments. Observe that the lemma below
actually gives a bit more information than what we need for Lemma
\ref{epsilon}. In what follows let $\Phi$ be the set of all
possible injective mappings $\phi: [h] \rightarrow [r]$, and
notice that there are $r!/(r-h)!$ elements in $\Phi$.

\begin{lemma}
\label{main} For every $h > 2$ there exists
$N_{\ref{main}}=N_{\ref{main}}(h)$ so that the following holds.
Let $H$ be a fixed graph with $m$ edges and $h$ vertices. If $r
\geq N_{\ref{main}}$ and $W$ is a labeled weighted graph on $r$
vertices satisfying $W(\phi)= \delta_H(p)$ for all $\phi \in
\Phi$, then $w(i,j) \in \{p,\overline{p}\}$ for all $1 \leq i < j
\leq r$. Furthermore, if $\gcd({h \choose 2},m)=1$ then all edge
weights are the same, and if $\gcd({h \choose 2},m)>1$ then either
all edge weights are the same, or else there exists one vertex
whose deletion from $W$ yields a subgraph with $r-1$ vertices all
of whose edge weights are the same.
\end{lemma}

We split the proof of Lemma \ref{main} into two parts. We
initially prove Lemma \ref{main1} below showing that all $w(i,j)$
are taken from $\{p,\overline{p}\}$. We then use this lemma in
order to show that in fact most of the $w(i,j)$ are either $p$ or
$\overline{p}$.

An important ingredient in the proof of Lemma \ref{main1} will be
a theorem of Gottlieb \cite{G}, concerning the rank of set
inclusion matrices. For integers $r \ge h > 2$, the {\em inclusion
matrix} $A(r,h)$ is defined as follows: The rows of $A(r,h)$ are
indexed by $h$-element subsets of $[r]$, and the columns by the
$2$-element subsets of $[r]$. Entry $(i,j)$ of $A(r,h)$ is $1$ if
the $2$-element set, whose index is $j$, is contained in the
$h$-element set, whose index is $i$. Otherwise, this entry is $0$.
Notice that $A(r,h)$ is a square matrix if and only if $r=h+2$,
and that for $r > h+2$, $A(r,h)$ has more rows than columns.
Trivially, $rank(A(r,h)) \leq {r \choose 2}$. However, Gottlieb
\cite{G} proved \footnote{Gottliebs's theorem actually deals with
the more general case where the columns are indexed by the $d$
element subsets of $[r]$ where $2 \leq d \leq h$, and in that case
the rank is ${r \choose d}$ for all $r \geq h+d$.} that in fact

\begin{theo}[Gottlieb \cite{G}]
\label{gottlieb} $rank(A(r,h)) = {r \choose 2}$ for all $r \ge
h+2$.
\end{theo}
\begin{lemma}
\label{main1}
Let $H$ be a fixed graph with $h > 2$ vertices. If $r \ge h+2$
and $W$ is a labeled weighted graph on $r$
vertices satisfying $W(\phi)= \delta_H(p)$ for all $\phi \in
\Phi$, then $w(i,j) \in \{p,\overline{p}\}$ for all $1 \leq i < j
\leq r$.
\end{lemma}

\paragraph{Proof:}
We associate a variable $x_{i,j}$ for each $1 \leq i < j \leq r$,
which represents the unknown $w(i,j)$. Thus, for any $\phi \in
\Phi$ we have that $W(\phi)$ is given by the polynomial
\begin{equation}\label{DefPphi}
P_\phi = \prod_{(i,j) \in E(H)}x_{\phi(i),\phi(j)} \prod_{(i,j)
\not \in E(H)}(1-x_{\phi(i),\phi(j)})\;.
\end{equation}
As our assumption is that $W(\phi)=\delta_H(p)$ for all $\phi \in
\Phi$ we have the following set of $r!/(r-h)!$ polynomial
equations $E_{\phi}$:
\begin{equation}\label{PolyEq}
E_{\phi}:~~~~ \prod_{(i,j) \in E(H)}x_{\phi(i),\phi(j)}
\prod_{(i,j) \not \in E(H)}(1-x_{\phi(i),\phi(j)})=\delta_H(p)\;.
\end{equation}
Our goal now is to show that the only solution to this system is
$x_{i,j} \in \{p,\overline{p}\}$.

For a vertex set $S \subseteq [r]$ of size $h$, let $E[S]$ denote
the ${h \choose 2}$ edges of $W$ induced by $S$. Let ${\cal S}$ be
the set of all $h$-element subsets of $V(W)=[r]$ and notice that
$|{\cal S}|={r \choose h}$. For every $S \in {\cal S}$ let
$\Phi_S$ be the set of all $h!$ elements of $\Phi$ that are
bijections on $S$. For every set $S$ let us take the product of
the $h!$ equations $\{E_{\phi}: \phi \in \Phi_S\}$ of
(\ref{PolyEq}). We thus get the following system of ${r \choose
h}$ polynomial equations (one for every $S \in {\cal S}$) with ${r
\choose 2}$ variables (one for each $(i,j) \in E(W)$):
\begin{equation}\label{DefPS1}
E_S:~~~~ \prod_{(i,j) \in E[S]} \left(x_{i,j}^m \cdot
(1-x_{i,j})^{{h \choose 2} -m}\right)^{2!(h-2)!} =
(\delta_H(p))^{h!}\;.
\end{equation}
In order to show that the only solution of the equations $E_S$ is
given by $x_{i,j} \in \{p,\overline{p}\}$, it would be convenient
to first transform them to linear equalities, by taking logarithms
on both sides. Define
\begin{equation}\label{Defy}
y_{i,j} = m \cdot \log (x_{i,j}) + \left({h \choose
2}-m\right)\cdot \log (1-x_{i,j})\;
\end{equation}
and note that if we take logarithm of the equations given in
(\ref{DefPS1}) and use the $y_{i,j}$ defined above, we thus obtain
an equivalent system of {\em linear} equations on the ${r \choose
2}$ variables $y_{i,j}$, where equation $E_S$ becomes
\begin{equation}\label{TheEquation}
E'_S:~~~~\sum_{(i,j) \in E[S]} y_{i,j} = {h \choose 2} \cdot
\log(\delta_H(p))\;.
\end{equation}
We can write the ${r \choose h}$ linear equations $E'_S$ as $Ax=b$
where $A$ in an ${r \choose h} \times {r \choose 2}$ matrix, and
$b$ is the all ${h \choose 2} \cdot \log(\delta_H(p))$ vector. A
key observation at this point is that $A$ is precisely the
inclusion matrix $A(r,h)$. Since $r \ge h+2$ we obtain, by Theorem
\ref{gottlieb}, that the system has a unique solution and the
values of the variables $y_{i,j}$ are uniquely determined. Now, as
each set $S \in {\cal S}$ is of size $h$ it is clear that setting
$y_{i,j} = \log (\delta_H(p))$ for all $(i,j)$ gives a valid
solution of the linear equations given in (\ref{TheEquation}), and
by the above observation, this is in fact the unique solution.
Recalling the definition of $y_{i,j}$ in (\ref{Defy}), this
implies that for all $(i,j)$ we have
$$
x_{i,j}^m \cdot (1-x_{i,j})^{{h \choose 2} - m}=\delta_H(p)\;.
$$
Now, as $\{p,\overline{p}\}$ are the only solutions to the above
equation, we deduce that indeed $x_{i,j} \in\{p,\overline{p}\}$,
proving the lemma. \qed

\bigskip

For the proof of Lemma \ref{main}, we will need another simple
lemma. A graph is called {\em pairwise regular} if there exists a
number $t$ so that $d(x)+d(y)-d(x,y)=t$ for all pairs of distinct
vertices $x,y$. Here $d(v)$ denotes the degree of $v$ and
$d(u,v)=1$ if $(u,v)$ is an edge, otherwise $d(u,v)=0$. A graph is
called {\em pairwise outer-regular} if there exists a number $t$
so that $d(x)+d(y)-2d(x,y)=t$ for all pairs of distinct vertices
$x,y$. Trivially, a graph is pairwise regular if and only if its
complement is. The same holds for pairwise outer-regular. It is
also trivial that the complete graph (and the empty graph) is both
pairwise regular and pairwise outer-regular. Notice, that
$K_{1,2}$ is also pairwise regular, and that $K_{1,3}$ is also
pairwise outer-regular. The following lemma, whose proof is
deferred to the end of this section, establishes that these are
the only non-trivial cases.

\begin{lemma}
\label{pairwise-regular} The only non-complete and non-empty
graphs which are pairwise regular are $K_{1,2}$ and its
complement. The only non-complete and non-empty graphs which are
pairwise outer-regular are $K_{1,3}$ and its complement.
\end{lemma}

\paragraph{Proof of Lemma \ref{main}:} Notice that if
$p=\overline{p}$ then there is actually nothing to prove, since
Lemma \ref{main1} already yields the desired conclusion. Hence,
assume $p \neq \overline{p}$. Observe, that this implies that $H$
in not the complete graph nor the empty graph as in these two
cases $p = \overline{p}$.

By Lemma \ref{main1}, each edge weight is either $p$ or
$\overline{p}$. We color the edges of $W$ with two colors: blue
for edges whose weight is $p$ and red for edges whose weight is
$\overline{p}$. We may assume that our coloring is non-trivial,
that is, that we have both red and blue edges, since otherwise
there is nothing to prove. Each $\phi \in \Phi$ defines a labeled
copy of $H$ in $W$. Let $b(\phi)$ be the number of edges of $H$
mapped to blue edges and let $a(\phi)$ be the number of
non-edges\footnote{The non-edges of $H$ are all the pairs $i,j$
that are not connected in $H$.} of $H$ mapped to blue edges. Then,
the number of edges of $H$ mapped to red edges is $m-b(\phi)$ and
the number of non-edges of $H$ mapped to red edges is ${h \choose
2}-m-a(\phi)$. Thus, we have for every $\phi \in \Phi$ that
\begin{equation}
\label{e1}
\delta_H(p)=W(\phi)=p^{b(\phi)}\overline{p}^{m-b(\phi)}(1-p)^{a(\phi)}(1-\overline{p})^{{h
\choose 2}-m-a(\phi)}\;.
\end{equation}
Multiplying (\ref{e1}) by
$\overline{p}^{b(\phi)}(1-\overline{p})^{a(\phi)}$ we get that
\begin{equation}\label{e0}
\delta_H(p) \cdot \overline{p}^{b(\phi)}(1-\overline{p})^{a(\phi)}
= p^{b(\phi)}(1-p)^{a(\phi)}\overline{p}^m(1-\overline{p})^{{h
\choose 2}-m} =\delta_H(p) \cdot p^{b(\phi)}(1-p)^{a(\phi)}\;
\end{equation}
where in the second equality we use the fact that
$\overline{p}^m(1-\overline{p})^{{h \choose 2}-m} =\delta_H(p)$.
This implies that
\begin{equation}
\label{e2}
\left(\frac{p}{\overline{p}}\right)^{b(\phi)}\left(\frac{1-p}{1-\overline{p}}\right)^{a(\phi)}=1\;.
\end{equation}
On the other hand, since $p$ and $\overline{p}$ are both solutions
of the equation $x^m \cdot (1-x)^{{h\choose 2}-m}=\delta_H(p)$ we
also know that
\begin{equation}
\label{e3}
\left(\frac{p}{\overline{p}}\right)^m\left(\frac{1-p}{1-\overline{p}}\right)^{{h
\choose 2}-m}=1\;.
\end{equation}
Thus, solving (\ref{e2}) for $\frac{p}{\overline{p}}$ and plugging
it into (\ref{e3}) gives that for any $\phi$
\begin{equation}
\label{e4} a(\phi)\cdot m = b(\phi)\cdot \left({h \choose
2}-m\right)\;.
\end{equation}

Consider first the case where $\gcd(m,{h \choose 2})=1$. This
implies that $\gcd(m,{h \choose 2}-m)=1$. Since the red-blue
coloring is not trivial there is a $K_h$ subgraph of $W$ which
contains both red and blue edges. Thus there exists $\phi \in
\Phi$ so that $0 < a(\phi) + b(\phi) < {h \choose 2}$. There are
two ways in which (\ref{e4}) can be satisfied: the first is if
$a(\phi)=b(\phi)=0$, but this violates the fact that $0 < a(\phi)
+ b(\phi)$. The second is if $a(\phi)$ is a multiple of ${h
\choose 2}-m$ and $b(\phi)$ is a multiple of $m$, but this
violates $a(\phi) + b(\phi) < {h \choose 2}$. Thus, the coloring
must be trivial, and we are done.

Now consider the case $\gcd(m,{h \choose 2})> 1$. By Ramsey's
Theorem if $N_{\ref{main}}(h)$ is sufficiently large, there is a
monochromatic copy of $K_{3h-8}$ in $W$. Let $T$ denote a maximal
monochromatic copy in $W$. Thus, $T$ has $t$ vertices and $r
> t \ge 3h-8$. Suppose, wlog, that $T$ is completely red. Let
$x$ be a vertex outside $T$. By maximality of $T$, there exists $y
\in T$ so that $(x,y)$ is blue. Suppose $x$ has at least $h-2$ red
neighbors in $T$, say $(x,v_1),\ldots,(x,v_{h-2})$ are all red.
Then, $\{x,y,v_1,\ldots,v_{h-2}\}$ induce a copy of $K_h$ which
has precisely one blue edge. If $\phi$ is any bijection onto this
copy then $a(\phi)+b(\phi)=1$, but this must violate (\ref{e4})
and hence the coloring must be trivial and we are done.

We may now assume that each vertex $x$ outside $T$ has at most
$h-3$ red neighbors in $T$. Now, if $t=r-1$ then there are at most
$r-1$ blue edges in our coloring, all incident with $x$, and we
are done. Otherwise, there are at least two vertices $x_1$ and
$x_2$ outside $T$, that have at least $t-2(h-3) \ge 3h-8-2h+6=h-2$
common neighbors $\{v_1,\ldots,v_{h-2}\}$ in $T$ so that all edges
$(x_i,v_j)$ are blue for $i=1,2$ and $j=1,\ldots,h-2$.

Consider first the case where $(x_1,x_2)$ is blue. Since
$\gcd(m,{h \choose 2})> 1$ we must have that $H$ is not $K_{1,2}$
nor its complement. Thus, by Lemma \ref{pairwise-regular}, $H$ is
not pairwise regular\footnote{Remember that at this point we know
that $H$ is neither a complete graph nor an edgeless graph.}. Let
$\{u_1,u_2\}$ and $\{u_3,u_4\}$ be two pairs of distinct vertices
of $H$ so that
\begin{equation}\label{violating}
d(u_1)+d(u_2)-d(u_1,u_2) \neq d(u_3)+d(u_4)-d(u_3,u_4)\;.
\end{equation}
Let $\phi_1$ be a bijection from $V(H)$ to
$\{x_1,x_2,v_1,\ldots,v_{h-2}\}$ mapping $u_1$ to $x_1$ and $u_2$
to $x_2$. Clearly, $b(\phi_1)= d(u_1)+d(u_2)-d(u_1,u_2)$.
Similarly, if $\phi_2$ is a bijection from $V(H)$ to
$\{x_1,x_2,v_1,\ldots,v_{h-2}\}$ mapping $u_3$ to $x_1$ and $u_4$
to $x_2$ then $b(\phi_2)= d(u_3)+d(u_4)-d(u_3,u_4)$. In
particular, we get from (\ref{violating}) that $b(\phi_1) \neq
b(\phi_2)$. We claim however that this is impossible as in fact
$b(\phi_1) = b(\phi_2)$. Indeed, by combining (\ref{e4}) for
$\phi_1$ and for $\phi_2$ we get that
$a(\phi_1)/a(\phi_2)=b(\phi_1)/b(\phi_2)$. Further we have
$b(\phi_1)+a(\phi_1)=b(\phi_2)+a(\phi_2)$ as both sides are equal
to the number of blue edges in the corresponding induced $K_h$ of
$W$. Combining the two equations we get $b(\phi_1) = b(\phi_2)$.

Consider finally the case where $(x_1,x_2)$ is red. Assume first
that $H$ is not $K_{1,3}$ nor its complement. Thus, by Lemma
\ref{pairwise-regular}, $H$ is not pairwise outer-regular. Let
$\{u_1,u_2\}$ and $\{u_3,u_4\}$ be two pairs of distinct vertices
of $H$ so that
$$
d(u_1)+d(u_2)-2d(u_1,u_2) \neq d(u_3)+d(u_4)-2d(u_3,u_4)\;.
$$
Let $\phi_1$ be a bijection from $V(H)$ to
$\{x_1,x_2,v_1,\ldots,v_{h-2}\}$ mapping $u_1$ to $x_1$ and $u_2$
to $x_2$. Clearly, $b(\phi_1)= d(u_1)+d(u_2)-2d(u_1,u_2)$.
Similarly, if $\phi_2$ is a bijection from $V(H)$ to
$\{x_1,x_2,v_1,\ldots,v_{h-2}\}$ mapping $u_3$ to $x_1$ and $u_4$
to $x_2$ then $b(\phi_2)= d(u_3)+d(u_4)-2d(u_3,u_4)$. In
particular, $b(\phi_1) \neq b(\phi_2)$. As in the previous case,
this is a contradiction. If $H=K_{1,3}$ then $h=4$ and we can use
the fact that $x_1$ has at least $3h-8-(h-3)=3$ blue neighbors in
$T$ denoted $y_1,y_2,y_3$. Thus, $x_1,y_1,y_2,y_3$ have a red
triangle and a blue $K_{1,3}$. Let $\phi_1$ map the vertex of
degree 3 of $H$ to $x_1$ and the rest to $y_1,y_2,y_3$, yielding
$b(\phi_1)=3$. Let $\phi_2$ map the vertex of degree 3 of $H$ to
$y_1$ and the rest to $x_1,y_2,y_3$ yielding $b(\phi_2)=1$. Again,
$b(\phi_1) \neq b(\phi_2)$, a contradiction. The case of the
complement of $K_{1,3}$ is proved in the same way. $\qed$

\bigskip

For the proof of Lemma \ref{epsilon}, we will need the following
simple fact
\begin{claim}\label{stable}
For any integer $p$ and $\delta$ there is a
$\gamma=\gamma_{\ref{stable}}(\delta,p)$ with the following
property: Let $A$ be any $p \times p$ non-singular $0/1$ matrix, let
$b$ be any vector in $\mathbb{R}^p$ and let $x \in \mathbb{R}^p$ be
the unique solution of the system of linear equations $Ax=b$. Then
if $b'$ satisfies $\ell_{\infty}(b',b) \leq \gamma$ then the unique
solution $x'$ of $Ax'=b'$ satisfies $\ell_{\infty}(x',x) \leq
\delta$.
\end{claim}

\paragraph{Proof:} Fix any $p \times p$ non-singular matrix $A$ with
$0/1$ entries. Then the solution of $Ax=b$ is given by $x=A^{-1}b$.
As $x_i=\sum_{j=1}^p A_{i,j}^{-1} \cdot b_j$ is a continuous
function of $b$ it is clear to for any $\delta$ there is a
$\gamma=\gamma(\delta,A)$ such that if $\ell_{\infty}(b',b) \leq
\gamma$ then the unique solution $x'$ of $Ax'=b'$ satisfies
$\ell_{\infty}(x',x) \leq \delta$. Now, as there are finitely many
$0/1$ $p \times p$ matrices, we can set
$\gamma=\gamma_{\ref{stable}}(\delta,p)=\min_A \gamma(\delta,A)$,
where the minimum is taken over all $0/1$ $p \times p$ matrices.
$\qed$

\paragraph{Proof of Lemma \ref{epsilon}:} The lemma is an
immediate consequence of Lemma \ref{main} using standard arguments
of continuity; the continuity of polynomials as functions, and the
continuity of unique solutions to linear systems that is given in
Lemma \ref{stable} above. First we can take
$N_{\ref{epsilon}}(h)=N_{\ref{main}}(h)$. Now, given any $r \geq
N_{\ref{main}}(h)$ and $\epsilon$ we need to show that if all
$W(\phi)$ are very close to $\delta_{H}(p)$ then we can get the
conclusion of Lemma \ref{epsilon}.

First, we see that in Lemma \ref{main1} if all $W(\phi)$ are close
to $\delta_{H}(p)$ then by Lemma \ref{stable} any solution to the
linear equations $E'_S$ given in (\ref{TheEquation}) satisfies that
all $y_{i,j}$ are very close to $\log(\delta_H(p))$. By continuity
of $2^x$ this means that $x^m_{i,j}(1-x_{i,j})^{{h \choose 2}-m}$ is
close to $\delta_H(p)$, which again by continuity of $x^k$ implies
that either $x_{i,j}$ is close to $p$ or to $\overline{p}$. So the
conclusion of Lemma \ref{main1} is that if all $W(\phi)$ are close
to $\delta_{H}(p)$, then all densities are indeed close to either
$p$ or $\overline{p}$.

For the rest of the proof, in equations (\ref{e1}) and (\ref{e0}) we
replace $p$ and $\overline{p}$ with quantities close to them. This
means that (\ref{e2}) and (\ref{e3}) are no longer equations but
approximately equal to 1. This implies that in (\ref{e4}) we also
have approximate equality. However, note that as both sides of
(\ref{e4}) involve integers, once the two sides are close enough,
they must in fact be equal. Now, as the rest of the proof only
relies on the validity of (\ref{e4}) it follows verbatim as in the
proof of Lemma \ref{main}. $\qed$

\bigskip

It is interesting to note that we cannot hope to prove a stronger
version of Lemma \ref{main} in which {\em all} edge weights are
the same, regardless of $\gcd({h \choose 2}, m)$. Indeed, consider
the case where $H=C_h$ is a cycle with $h \ge 4$ vertices. For
every $r \ge h+1$, there are weighted complete graphs $W$ with $r$
vertices having $W(\phi)=\delta_H(p)$ for each $\phi \in \Phi$,
while still some edges of $W$ have weight $p$ and others have
weight $\overline{p}$. Indeed, assume that all weights of edges
not incident with $r \in W$ have weight $p$, and the $r-1$ edges
incident with $r$ have weight $\overline{p}$. Now, if the image of
$\phi$ does not contain $r$ then, clearly,
$$
W(\phi)=p^h(1-p)^{{h \choose 2}-h}=\delta_H(p).
$$
On the other hand, if the image of $\phi$ contains $r$ then
$$
W(\phi)=p^{h-2}\overline{p}^2(1-p)^{{h \choose 2}-2h+3}(1-\overline{p})^{h-3}.
$$
But note that $\overline{p}^2(1-\overline{p})^{h-3}$ is just
$\delta_H(p)^{2/h}$, and hence it also equals $p^2(1-p)^{h-3}$.
Consequently, $W(\phi)=\delta_H(p)$ in this case as well.

\paragraph{Proof of Lemma \ref{pairwise-regular}:} Let us say that
$(x_1,y_1)(x_2,y_2)$ are violating, with respect to the property of
being pairwise regular if $d(x_1)+d(y_1)-d(x_1,y_1) \neq
d(x_2)+d(y_2)-d(x_2,y_2)$ and violating with respect to the property
of being pairwise outer-regular if $d(x_1)+d(y_1)-2d(x_1,y_1) \neq
d(x_2)+d(y_2)-2d(x_2,y_2)$. Suppose first that $G$ is a pairwise
regular graph which is neither complete nor empty. We claim that
this implies that $|d(x)-d(y)| \le 1$ for any two vertices $x,y \in
V(G)$. Indeed, if there is a pair that violates this, then
$(x,z),(y,z)$ is violating for any $z$. Note that $G$ cannot be
regular, otherwise $(x_1,y_1)(x_2,y_2)$ is violating whenever
$(x_1,y_1)$ is and edge and $(x_2,y_2)$ is not. So partition the
vertices of $G$ into two non empty sets, $V_1$ and $V_2$, where all
the vertices of $V_1$ have degree $s$ and those of $V_2$ have degree
$s-1$.

If $|V_1| > 1$ then $V_1$ must be a clique otherwise
$(x_1,x_2)(x_2,y)$ is violating for any non-adjacent $x_1, x_2 \in V_1$ and $y \in
V_2$. In particular, we have $t=2s-1$. We also have that $|V_2|=1$
as otherwise $(x_1,x_2)(y_1,y_2)$ is violating for any $x_1, x_2 \in
V_1$ and $y_1,y_2 \in V_2$. If the unique vertex $v$ of $V_2$ is
connected to $x_1 \in V_1$ but not to $x_2 \in V_1$ then
$(x_1,v)(x_2,v)$ is clearly violating, so $v$ is either connected to
all the vertices of $V_1$ or else is an isolated vertex. If $v$ is
an isolated vertex then $s=1$, which implies that $V_1$ is a clique
of size 2, and $G$ is thus the complement of $K_{1,2}$. If $v$ is
connected to all the vertices of $V_1$ then $|V_1|=s-1$ which is
impossible, since in a graph with $s$ vertices there cannot be
vertices with degree $s$. If $|V_1| =1$ then we must have $|V_2| >
1$. Note that in this case $V_2$ must span an independent set as
otherwise $(x,y_1)(y_1,y_2)$ is violating for any choice of adjacent $y_1,y_2
\in V_2$ and $x \in V_1$. As $G$ is not edgeless we infer that
$s-1=1$ implying that $G$ is $K_{1,2}$.

Suppose now that $G$ is pairwise outer-regular and is neither
complete nor empty. Following the same reasoning as above, we must
have for any two vertices $x,y$ of $G$, that $|d(x)-d(y)| \le 2$.
Again, note that $G$ cannot be regular, so partition the vertices of
$G$ into two non-empty sets, $V_1$ and $V_2$, where all of the
vertices of $V_1$ have degree $s$ and all the vertices of $V_2$ have
degree $s-1$ or $s-2$. If $|V_1| > 1$ then again $V_1$ must span a
clique, as otherwise $(x_1,x_2)(x_2,y)$ is violating for any $x_1,x_2
\in V_1$ and $y \in V_2$, and therefore $t=2s-2$. Note that if $y
\in V_2$ is connected to $x \in V_1$ then $(x_1,x)(x,y)$ is
violating for any other $x_1 \in V_1$. Also, if $|V_2| \geq 2$ then
any pair of vertices of $V_2$ must be disconnected with degree $s-1$
as otherwise $(x_1,x_2)(y_1,y_2)$ is violating for any $x_1,x_2 \in
V_1$ and $y_1,y_2 \in V_2$. We thus get that the degree of vertices
in $V_2$ is zero, hence either $s=1$ or $s=2$. The former case
implies that $t=0$ and that $|V_1|=2$. This means that $G$ has just
one edge, which is not an outer-regular graph. If $s=2$ then
$|V_1|=3$ implying that $G$ is a triangle plus some isolated
vertices. If there is one such vertex then $G$ is the complement of
$K_{1,3}$, and if there are two such vertices $y_1,y_2$, then
$(y_1,y_2)(y_1,v)$ is violating for any $v \in V_1$. So assume that
$|V_1|=1$, which implies that $|V_2| \geq 2$. Let $x$ be the unique
vertex of $V_1$, and observe that if $x$ is connected to $y_1 \in
V_2$ but not to $y_2 \in V_2$ then $(x,y_1)(x,y_2)$ is violating. So
either $v$ is connected to all the vertices of $V_2$ or to none of
them, but note that the latter case is impossible as $s > s-1 \geq
0$. We now claim that $V_2$ must be edgeless. Indeed if $y_1,y_2 \in
V_2$ are connected and $d(y_1) \geq d(y_2)$ then $(y_1,y_2)(x,y_1)$
is violating. We infer that the degree of the vertices of $V_2$ is
1, so $G$ is either $K_{1,2}$, which is not outer regular, or
$K_{1,3}$. $\qed$

\end{document}